# FUNCTIONALS OF DIRICHLET PROCESSES, THE CIFARELLI–REGAZZINI IDENTITY AND BETA-GAMMA PROCESSES[1]

By Lancelot F. James

*Hong Kong University of Science and Technology*

Suppose that $P_\theta(g)$ is a linear functional of a Dirichlet process with shape $\theta H$, where $\theta > 0$ is the total mass and $H$ is a fixed probability measure. This paper describes how one can use the well-known Bayesian prior to posterior analysis of the Dirichlet process, and a posterior calculus for Gamma processes to ascertain properties of linear functionals of Dirichlet processes. In particular, in conjunction with a Gamma identity, we show easily that a generalized Cauchy–Stieltjes transform of a linear functional of a Dirichlet process is equivalent to the Laplace functional of a class of, what we define as, Beta-Gamma processes. This represents a generalization of an identity due to Cifarelli and Regazzini, which is also known as the Markov–Krein identity for mean functionals of Dirichlet processes. These results also provide new explanations and interpretations of results in the literature. The identities are analogues to quite useful identities for Beta and Gamma random variables. We give a result which can be used to ascertain specifications on $H$ such that the Dirichlet functional is Beta distributed. This avoids the need for an inversion formula for these cases and points to the special nature of the Dirichlet process, and indeed the functional Beta-Gamma calculus developed in this paper.

**1. Introduction.** Let $P$ denote a Dirichlet random probability measure on a Polish space $\mathcal{Y}$, with law denoted as $\mathcal{D}(dP|\theta H)$, where $\theta$ is a nonnegative scalar and $H$ is a (fixed) probability measure on $\mathcal{Y}$. In addition, let $\mathcal{M}$ denote the space of boundedly finite measures on $\mathcal{Y}$. This space contains the space of probability measures on $\mathcal{Y}$. The Dirichlet process was first made popular in Bayesian nonparametrics by Ferguson (1973) [see also Freedman

Received October 2003; revised July 2004.
[1]Supported in part by RGC Grant HKUST-6159/02P of the HKSAR.

([1963](#)) for an early treatment], and has subsequently been used in numerous statistical applications. Additionally, the Dirichlet process arises in a variety of interesting contexts outside of statistics. Formally, $P$ is said to be a Dirichlet process if and only if for each finite collection of disjoint measurable sets $A_1, \ldots, A_k$, the random vector $P(A_1), \ldots, P(A_k)$ has a Dirichlet distribution with parameters $\theta H(A_1), \ldots, \theta H(A_k)$. In particular, $P(A)$ is a Beta random variable for any measurable set $A$. An important representation of the Dirichlet, which is analogous to Lukacs characterization of the Gamma distribution, is $P(\cdot) = \mu(\cdot)/T$ where $\mu$ is a Gamma process with finite shape parameter $\theta H$ and $T = \int_{\mathcal{Y}} \mu(dy)$ is a Gamma random variable with shape $\theta$ and scale 1. The law of the Gamma process is denoted as $\mathcal{G}(d\mu|\theta H)$ and is characterized by its Laplace functional

$$\int_{\mathcal{M}} e^{-\mu(g)} \mathcal{G}(d\mu|\theta H) = e^{-\int_{\mathcal{Y}} \log[1+g(y)]\theta H(dy)}$$

for each positive bounded measurable function $g$ on $\mathcal{Y}$. For our purposes we shall consider the more general class of real-valued functions $g$ which satisfy the constraint

$$\int_{\mathcal{Y}} \log[1+|g(y)|]\theta H(dy) < \infty. \tag{1}$$

This condition, (1), as shown by Doss and Sellke (1982) and Feigin and Tweedie (1989), is necessary and sufficient for the existence of the linear functionals

$$P(g) = \int_{\mathcal{Y}} g(y) P(dy).$$

An important fact is that $T$ and $P$ are independent, which as we shall see, has a variety of implications.

An interesting problem initiated in a series of papers by Cifarelli and Regazzini (1990) is the study of the exact distribution of linear functionals $P(g)$ of the Dirichlet process. One of their contributions is the important identity

$$\int_{\mathcal{M}} \frac{1}{(1+zP(g))^\theta} \mathcal{D}(dP|\theta H) = e^{-\int_{\mathcal{Y}} \log[1+zg(y)]\theta H(dy)}, \tag{2}$$

where typically $z$ is in the complex plane $\mathscr{C}$. We call (2) the *Cifarelli–Regazzini identity*. The result in (2) is sometimes called the *Markov–Krein identity for means of Dirichlet processes*. Diaconis and Kemperman (1996) discuss some consequences of this result, which has implications relative to the Markov moment problem, continued fractions theory, exponential representations of analytic functions, and so on [see Kerov (1998) and Vershik, Yor and Tsilevich (2001)]. Vershik, Yor and Tsilevich (2001) expand upon this, emphasizing that the right-hand side of (2) is the Laplace functional of a Gamma process with shape $\theta H$. That is,

$$\int_{\mathcal{M}} \frac{1}{(1+zP(g))^\theta} \mathcal{D}(dP|\theta H) = \int_{\mathcal{M}} e^{-z\mu(g)} \mathcal{G}(d\mu|\theta H). \tag{3}$$

Their interpretation, which is in the sense of a Markov–Krein transform, is that the generalized Cauchy–Stieltjes transform of order $\theta$ of $P(g)$, where $P$ is a Dirichlet process with shape $\theta H$, is the Laplace transform of $\mu(g)$ when $\mu$ is the Gamma process with shape $\theta H$. The authors then exploit this fact to rederive (3) via an elementary proof using the independence property of $P$ and $T$. An interesting question is, what can one say about

$$(4) \qquad \int_{\mathcal{M}} \frac{1}{(1+zP(g))^q} \mathcal{D}(dP|\theta H)$$

when $\theta$ and $q$ are arbitrary positive numbers? That is, can one establish a relationship of (4) to the Laplace functional of some random measure, say $\mu^*$, which is similar to $\mu$, for all $q$ and $\theta$? Lijoi and Regazzini (2004) establish analytic results for (4), relating them to the Lauricella theory of multiple hypergeometric functions. Theorem 5.2 of their work gives analogues of (2), stating what they call a *Lauricella identity*, but does not specifically state a relationship such as (3). We should say for the case $\theta > q$ that it would not be terribly difficult to deduce an analogue of (4) from their result. However, this is not the case when $\theta < q$, which is expressed in terms of contour integrals. Their representations, for all $\theta$ and $q$, as clearly demonstrated by the authors, however have practical utility in regards to formulae for the density of $P(g)$. In this case, one wants to have an expression for (4), when $q = 1$ and for all $\theta$. Some related works include the papers of Kerov and Tsilevich (1998), Regazzini, Guglielmi and Di Nunno (2002), Regazzini, Lijoi and Prünster (2003) and the manuscript of James (2002).

1.1. *Preliminaries and outline.* In this paper we develop results that are complementary to the work of Lijoi and Regazzini (2004) and Vershik, Yor and Tsilevich (2001). In particular, we show that (3), as interpreted in Vershik, Yor and Tsilevich (2001), extends to a relationship between (4) and the Laplace functional of a class of what we call Beta-Gamma processes defined by scaling the Gamma process law by $T^{-d}$, for all numbers $d$ such that $\theta - d > 0$, that is, processes with laws equal to

$$(5) \qquad \mathcal{BG}(d\mu|\theta H, d) = \frac{\Gamma(\theta)}{\Gamma(\theta - d)} T^{-d} \mathcal{G}(d\mu|\theta H).$$

In particular, our main result concerns the choice of $d = \theta - q$ for arbitrary positive numbers which are not necessarily equal. The approach relies in part on, in this case, mostly familiar Bayesian prior posterior calculus for Dirichlet and Gamma processes in conjunction with the usage of the following well-known *Gamma identity* for $q > 0$:

$$(6) \qquad T^{-q} = \frac{1}{\Gamma(q)} \int_0^\infty v^{q-1} e^{-vT} \, dv.$$

That is to say, purely analytic arguments are replaced by Bayesian arguments using the familiar results in Ferguson (1973), Lo (1984) and Antoniak (1974), thus giving the derivations a much more interpretable Bayesian flavor. More specifically, albeit less well known, we use the results in Lo and Weng (1989) as demonstrated for more general processes in James (2002). This bypasses the need, for instance, to verify certain integrability conditions and the usage of limiting arguments. Moreover, somewhat conversely to Lijoi and Regazzini (2004), we show how properties of the Dirichlet and Beta-Gamma processes yield easily interesting identities related to Lauricella and Liouville integrals [see Lijoi and Regazzini (2004) and Gupta and Richards (2001)]. Although we exploit the independence property of $T$ and $P$ to prove our results, our approach is quite different from the methods used by Vershik, Yor and Tsilevich (2001) to prove (3). While their proof is certainly elegant, it does not seem possible to extend to other processes. Our methods, however, are influenced by their proof of an analogous result for the two-parameter extension of the Dirichlet process [see Pitman (1996)] which relies on (6) and the fact that such processes are based on scaled laws. That is to say, we present an approach which is extendable to other models [see James (2002), Section 6]. However, for the Beta-Gamma processes defined in (5), the independence property between $T$ and $P$ translates into the property

$$(7) \quad \int_{\mathcal{M}} h(P)\mathcal{D}(dP|\theta H) = \int_{\mathcal{M}} h(P)\mathcal{G}(d\mu|\theta H) = \int_{\mathcal{M}} h(P)\mathcal{BG}(d\mu|\theta H, d)$$

for all integrable $h$. The property (7) seems to suggest that the Beta-Gamma process may not have much utility relative to calculations involving $P$; however, it is precisely this property that we shall exploit. In the next section we shall first develop, rather quickly, two supporting results concerning the calculus of Dirichlet and Beta-Gamma processes. We will then show how these results are used to easily derive our main results in Theorems 2.1 and 2.2 based on Bayesian arguments. We close the paper by showing how our methodology, a Beta-Gamma calculus for Dirichlet processes, leads to a functional analogue of the classical Beta-Gamma calculus for random variables. That is, we provide conditions on $H$ such that $P_\theta(g)$ is Beta distributed.

**2. Functionals of Dirichlet processes, the Cifarelli–Regazzini identity and Beta-Gamma processes.** We start by recalling some properties of the Dirichlet process. Let $Y_1, \ldots, Y_n$ denote random elements in the space $\mathcal{Y}$, which conditional on $P$ are i.i.d. with law $P$. $P$ is a Dirichlet process with shape $\theta H$. These specifications define a joint law of $(\mathbf{Y}, P)$, where, from Ferguson (1973), it follows that the posterior distribution of $P|\mathbf{Y}$ is also Dirichlet with shape

$$(\theta + n)H_n = \theta H + \sum_{i=1}^{n} \delta_{Y_i}.$$

Additionally, the marginal distribution of $\mathbf{Y}$ is the Blackwell and MacQueen (1973) distribution described as

$$\mathbb{P}(d\mathbf{Y}|\theta H) = \frac{\Gamma(\theta)}{\Gamma(\theta+n)}\theta H(dY_1)\prod_{i=2}^{n}\left(\theta H + \sum_{j=1}^{i-1}\delta_{Y_j}\right)(dY_i).$$

The Blackwell–MacQueen distribution admits ties. Hence one can represent $\mathbf{Y} = (\mathbf{Y}^*, \mathbf{p})$, where $\mathbf{Y}^* = \{Y_1^*, \ldots, Y_{n(\mathbf{p})}^*\}$ denotes the $n(\mathbf{p}) \leq n$ unique values of $\mathbf{Y}$. The expression $\mathbf{p} = \{C_1, \ldots, C_{n(\mathbf{p})}\}$ denotes a partition of the integers $\{1, \ldots, n\}$ with $n(\mathbf{p})$ elements. The $C_j = \{i : Y_i = Y_j^*\}$ for $j = 1, \ldots, n(\mathbf{p})$ denote the collection of values equal to each unique $Y_j^*$, for $j = 1, \ldots, n(\mathbf{p})$. The cardinality of each cell $C_j$ is denoted as $e_{j,n}$. When $H$ is nonatomic, then one can write

$$\mathbb{P}(dY|\theta H) = \pi(\mathbf{p}|\theta)\prod_{j=1}^{n(\mathbf{p})}H(dY_j^*),$$

where

$$\pi(\mathbf{p}|\theta) = \frac{\theta^{n(\mathbf{p})}\Gamma(\theta)}{\Gamma(\theta+n)}\prod_{j=1}^{n(\mathbf{p})}(e_{j,n}-1)!$$

is a variant of Ewens' (1972) sampling formula, which was independently derived by Antoniak (1974). It is also called the Chinese Restaurant process [see Pitman (1996)] and plays a fundamental role in Lo (1984). If $H$ is discrete with probability mass function $\rho$, then

$$P(d\mathbf{Y}|\theta H) = \frac{\Gamma(\theta)}{\Gamma(\theta+n)}\prod_{j=1}^{n(\mathbf{p})}\frac{\Gamma(\theta\rho(Y_j^*)+e_{j,n})}{\Gamma(\theta\rho(Y_j^*))}.$$

In any case, note that appealing to standard Bayesian arguments, the results of Ferguson (1973) imply that one has

(8) $$\int_{\mathcal{M}}h(P)\mathcal{D}(dP|\theta H) = \int_{\mathcal{Y}^n}\left[\int_{\mathcal{M}}h(P)\mathcal{D}(dP|(\theta+n)H_n)\right]\mathbb{P}(d\mathbf{Y}|\theta H).$$

This simple consequence is fundamental to our presentation. It is evident that (8) along with the various forms of $P(d\mathbf{Y}|\theta H)$ yield nontrivial expressions which might otherwise require an appeal to, for instance, the theory of special functions or combinatorics. In the same spirit, we now derive the posterior distribution of the Beta-Gamma processes. From Lo and Weng (1989), one has the following disintegration of measures:

(9) $$\prod_{i=1}^{n}\mu(dY_i)\mathcal{G}(d\mu|\theta H) = \frac{\Gamma(\theta+n)}{\Gamma(\theta)}\mathcal{G}\left(d\mu|\theta H + \sum_{i=1}^{n}\delta_{Y_i}\right)\mathbb{P}(d\mathbf{Y}|\theta H),$$

where $\mathcal{G}(d\mu|\theta H + \sum_{i=1}^{n} \delta_{Y_i})$ denotes a Gamma process with shape $\theta H + \sum_{i=1}^{n} \delta_{Y_i}$. Using (7), it is easy to see that $\int_{\mathcal{M}} \prod_{i=1}^{n} P(dY_i)\mathcal{BG}(d\mu|\theta H, d) = \mathbb{P}(d\mathbf{Y}|\theta H)$. Combining this fact with (9) easily yields the following description of the posterior distribution of a Beta-Gamma process.

PROPOSITION 2.1. *Let $\mu$ have law $\mathcal{BG}(d\mu|\theta H, d)$ defined for all $d$, such that $\theta - d > 0$. Then from (7), the law of $P = \mu/T$ is $\mathcal{D}(dP|\theta H)$. Suppose that $Y_1, \ldots, Y_n|P$ are i.i.d. $P$; then the posterior distribution of $\mu|\mathbf{Y}$ is a Beta-Gamma process with parameters $(\theta + n)H_n$ and $n + d$, defined as*

$$\mathcal{BG}(d\mu|(\theta+n)H_n, n+d) = \frac{\Gamma(\theta+n)}{\Gamma(\theta-d)} T^{-(n+d)} \mathcal{G}(d\mu|(\theta+n)H_n).$$

*Hence, similar to* (8), *one has*

$$
\begin{aligned}
(10) \quad & \int_{\mathcal{M}} f(\mu) \mathcal{BG}(d\mu|\theta H, d) \\
& = \int_{\mathcal{Y}^n} \left[ \int_{\mathcal{M}} f(\mu) \mathcal{BG}(d\mu|(\theta+n)H_n, n+d) \right] \mathbb{P}(d\mathbf{Y}|\theta H)
\end{aligned}
$$

*for all integrable $f$. Note that setting $d = 0$ shows that if $\mu$ is $\mathcal{G}(d\mu|\theta H)$, then its posterior distribution is $\mathcal{BG}(d\mu|(\theta+n)H_n, n)$, which is not a Gamma process.*

REMARK 2.1. Note that the use of (7), (8) and (10) sets up a myriad of interesting equivalences which will prove useful in our derivations. However, we do point out that while (10) implies (8), the converse is not true.

Another important property of the Gamma process that we shall exploit is the algebraic identity

$$(11) \int_{\mathcal{M}} e^{-(vT+w\mu(g))} \mathcal{G}(d\mu|\theta H) = (1+v)^{-\theta} e^{-\int_{\mathcal{Y}} \log[1+(w/(1+v))g(y)]\theta H(dy)}.$$

Let

$$\mathcal{B}(du|a,b) = \frac{\Gamma(a+b)}{\Gamma(a)\Gamma(b)} u^{a-1}(1-u)^{b-1} du \quad \text{for } 0 < u < 1$$

denote the density of a Beta random variable with parameters $(a, b)$. We now establish our final preliminary result before our main theorem.

PROPOSITION 2.2. *Let $\theta$ and $q$ be arbitrary nonnegative numbers. Then for any integer $n \geq 0$ that satisfies the constraint $\theta + n - q > 0$, the following formula holds:*

$$
\begin{aligned}
(12) \quad & \frac{\Gamma(\theta+n)}{\Gamma(\theta+n-q)} \int_{\mathcal{M}} \frac{1}{(T+z\mu(g))^q} \mathcal{G}(d\mu|(\theta+n)H_n) \\
& = \int_{\mathcal{M}} e^{-z\mu(g)} \mathcal{BG}(d\mu|(\theta+n)H_n, \theta+n-q).
\end{aligned}
$$

PROOF. Apply the Gamma identity to $(T+\mu(g))^{-q}$ and then (11) with $v=w$ and $\theta H$ replaced by $(\theta+n)H_n$, to show that the left-hand side of (12) is equal to

$$\frac{\Gamma(\theta+n)}{\Gamma(\theta+n-q)\Gamma(q)} \int_0^\infty v^{q-1} e^{-\int_{\mathcal{Y}} \log[1+(v/(1+v))zg(y)]\theta H(dy)}$$

$$\times \prod_{j=1}^{n(\mathbf{p})} \left(1 + \frac{v}{(1+v)}zg(Y_j^*)\right)^{-e_{j,n}} (1+v)^{-(\theta+n)} dv.$$

Similarly, the following expression is obtained for the right-hand side of (12) by applying the Gamma identity to $T^{-(\theta+n-q)}$, (11) with $w=1$ and a further change of variable:

$$(13) \quad \int_0^1 e^{-\int_{\mathcal{Y}} \log[1+uzg(y)]\theta H(dy)} \prod_{j=1}^{n(\mathbf{p})} (1+uzg(Y_j^*))^{-e_{j,n}} \mathcal{B}(du|q,\theta+n-q).$$

The result is obtained by applying the transformation $u=v/(1+v)$. □

We now present a new result which relates the generalized Cauchy–Stieltjes transform of Dirichlet process linear functionals to the Laplace functional of Beta-Gamma processes. This presents a generalization of the Cifarelli–Regazzini identity, complementary to the Lauricella identities deduced in Lijoi and Regazzini [(2004), Theorem 5.2]. We also present some interesting additional identities.

THEOREM 2.1. *Let $\mathcal{D}(dP|\theta H)$ denote a Dirichlet process with shape $\theta H$. Let $g$ denote a function satisfying (1). Then the following relationships are established:*

(i) *For any positive $q$ and $\theta$,*

$$(14) \quad \int_{\mathcal{M}} (1+zP(g))^{-q} \mathcal{D}(dP|\theta H) = \int_{\mathcal{M}} e^{-z\mu(g)} \mathcal{BG}(d\mu|\theta H, \theta-q).$$

*Note that the law $\mathcal{BG}(d\mu|\theta H, \theta-q)$ exists for all positive $\theta$ and $q$, and arbitrary $H$, since $\theta - (\theta-q) = q > 0$.*

(ii) *For any positive $q$ and $\theta$, and any integer $n \geq 0$ which satisfies $\theta + n - q > 0$, the quantities in (14) are equivalent to $\int_{\mathcal{Y}^n} [\int_{\mathcal{M}} e^{-\mu(g)} \mathcal{BG}(d\mu|(\theta+n)H_n, \theta+n-q)] \mathbb{P}(d\mathbf{Y}|\theta H)$. An explicit expression can be deduced from the equivalence of the inner term to (13). In particular, when $H$ is nonatomic, the expression is equivalent to*

$$\sum_{\mathbf{p}} \pi(\mathbf{p}|\theta) \int_0^1 e^{-\int_{\mathcal{Y}} \log[1+uzg(y)]\theta H(dy)}$$

$$\times \left[\prod_{j=1}^{n(\mathbf{p})} \int_{\mathcal{Y}} (1+uzg(y))^{-e_{j,n}} H(dy)\right] \mathcal{B}(du|q,\theta+n-q).$$

*For the Gamma process with law $\mathcal{G}(d\mu|\theta H)$, its Laplace functional may be represented as above for all $n \geq 1$ and $q = \theta$.*

(iii) *When $\theta - q > 0$, statements* (i) *and* (ii) *with $n = 0$ imply that*

$$\int_{\mathcal{M}} (1 + zP(g))^{-q} \mathcal{D}(dP|\theta H) \tag{15}$$
$$= \int_0^1 e^{-\int_{\mathcal{Y}} \log[1 + uzg(y)]\theta H(dy)} \mathcal{B}(du|q, \theta - q),$$

*which coincides with the result in Lijoi and Regazzini* [(2004), *Theorem* 5, *equation* (5.2)].

PROOF. A general strategy is formed by first writing $(1 + zP(g))^{-q} = T^q(T + z\mu(g))^{-q} = h(P)$. For the proof of statement (i), we first assume without loss of generality that $q = n + d$, where $d$ is a positive number such that $\theta - d > 0$, and $n \geq 0$ is an integer chosen such that $\theta + n - q > 0$. This means that $T^q = T^{n+d}$. Now using (8) and then (7) with $\mathcal{BG}(d\mu|(\theta+n)H_n, q)$ yields

$$\int_{\mathcal{M}} (1 + zP(g))^{-q} \mathcal{D}(dP|\theta H)$$
$$= \frac{\Gamma(\theta + n)}{\Gamma(\theta - d)} \int_{\mathcal{Y}^n} \left[ \int_{\mathcal{M}} \frac{1}{(T + z\mu(g))^q} \mathcal{G}(d\mu|(\theta + n)H_n) \right] \mathbb{P}(d\mathbf{Y}|\theta H).$$

Apply Proposition 2.2 to the inner term, recalling that $\theta + n - q = \theta - d$. This yields the desired expression,

$$\int_{\mathcal{Y}^n} \left[ \int_{\mathcal{M}} e^{-z\mu(g)} \mathcal{BG}(d\mu|(\theta + n)H_n, \theta + n - q) \right] \mathbb{P}(d\mathbf{Y}|\theta H) \tag{16}$$
$$= \int_{\mathcal{M}} e^{-z\mu(g)} \mathcal{BG}(d\mu|\theta H, \theta - q).$$

Note how again an appeal to a Bayesian argument, that is, using (10) in Proposition 2.1 with $f(\mu) = e^{-z\mu(g)}$, is used to deduce easily the equivalence of the right- and left-hand sides of (16). □

REMARK 2.2. Since $H$ is an arbitrary distribution, the result applies to a Dirichlet process posterior distribution based on, say, a sample of size $m$ having no particular relationship to $n$. For concreteness, suppose $P$ is a Dirichlet process with shape $\alpha P_0 + \sum_{i=1}^m \delta_{X_i}$, where $P_0$ is an arbitrary probability measure, $\alpha$ is a positive scalar and $X_1, \ldots, X_m$ are fixed points. Then the results in Theorem 2.1 hold for this $P$ by setting $\theta = \alpha + m$ and $\theta H = \alpha P_0 + \sum_{i=1}^m \delta_{X_i}$.

REMARK 2.3. As in Kerov and Tsilevich (1998) and Vershik, Yor and Tsilevich (2001), Theorem 2.1 applies to the joint distribution of linear functionals, say $(P(g_1), \ldots, P(g_m))$, where $g_1, \ldots, g_m$ are functions satisfying (1). The generalized Cauchy–Stieltjes transform for joint distributions is defined by replacing $zP(g)$ by $\sum_{i=1}^m z_i P(g_i)$. Since $\sum_{i=1}^m z_i P(g_i) = P(\sum_{i=1}^m z_i g_i)$, the result is seen by replacing $zg$ with $\sum_{i=1}^m z_i g_i$ in Theorem 2.1.

We now discuss some interesting results obtained from Theorem 2.1. Note the relative ease by which Bayesian arguments can be used to derive otherwise complex expressions such as that appearing in Theorem 2.1(ii). The case $q = 1$ is of particular interest in terms of giving an expression for the Cauchy–Stieltjes transform $P(g)$, which can be inverted to obtain an expression for the distribution of $P(g)$. Setting $q = 1$ in Theorem 2.1(ii) gives a variety of equivalent expressions which hold for all $\theta$ and $n \geq 1$. Here, as a corollary, we present the simplest expression that holds for all $\theta$ with $n = 1$.

COROLLARY 2.1. *Let $\mathcal{D}(dP|\theta H)$ denote a Dirichlet process with shape $\theta H$, where $H$ is an arbitrary probability measure. Let $g$ denote a function satisfying* (1); *then for all $\theta > 0$,*

$$
\begin{aligned}
(17) \quad &\int_{\mathcal{M}} \frac{\mathcal{D}(dP|\theta H)}{(1 + zP(g))} \\
&= \int_0^1 e^{-\int_{\mathcal{Y}} \log[1+uzg(y)]\theta H(dy)} \int_{\mathcal{Y}} \frac{\theta H(dy)}{1 + uzg(y)} (1-u)^{\theta-1} \, du.
\end{aligned}
$$

*When $\theta > 1$, this expression equates to the expression in Theorem* 2.1(iii). *When $\theta = 1$, Theorem* 2.1(ii) *shows that the right-hand side of* (17) *is the Laplace functional of a Gamma process with shape $H$, which corresponds to* (2).

The expression (17) can be seen as complementary to the expressions obtained in Lijoi and Regazzini (2004). However, our results are quite different in the case where $0 < \theta < 1$, where those authors obtained an expression in terms of contour integrals.

Let $\mathcal{L}(Z)$ denote the law of a random element $Z$. For the remainder of this work, let $\mu_{\theta,\theta-q}$ be a Beta-Gamma process with parameters $(\theta H, \theta - q)$, and let $U_{a,b}$ denote a Beta$(a,b)$ random variable. Let $T_\alpha$ denote a Gamma random variable with shape $\alpha$ and scale 1, and let $Y_1$ be a random element with distribution $H$. Let $\mu_\theta$ denote a Gamma process with shape $\theta H$ and assume that the variables $\mu_\theta, U_{a,b}, T_\alpha, Y_1$ are independent. Additionally, let $P_\theta$ denote a Dirichlet process with shape $\theta$. When convenient we will simply write $X = T$ to denote that the distribution of $X$ is equivalent to that of $T$. That is, $X = U_{a,b}$ means that $X$ has a Beta distribution with parameters $(a, b)$. The next result involves a series of distributional identities. These are based on Bayesian mixture representations deduced from the form of the posterior distribution mixed over the marginal distribution, $P(d\mathbf{Y}|\theta H)$. Some important consequences will be demonstrated thereafter.

THEOREM 2.2. *Let $\mu_{\theta,\theta-q}$ be a Beta-Gamma process with parameters $(\theta H, \theta - q)$ and let $\mu_\theta$ denote a Gamma process with shape $\theta H$. Then for all positive $\theta$ and $q$ and an integer $n$ chosen such that $\theta + n - q > 0$, the following distributional equalities hold:*

(i) *For all $\theta > 0$ and $q$, and an integer $n$ chosen such that $\theta + n - q > 0$,*

$$(18) \quad \mathcal{L}(\mu_{\theta,\theta-q}) = \mathcal{L}\left(U_{q,\theta+n-q}\mu_\theta + U_{q,\theta+n-q}\sum_{j=1}^{n(\mathbf{p})} G_{j,n}\delta_{Y_j^*}\right),$$

*where conditional on $\mathbf{p}$ the distinct variables on the right-hand side are mutually independent such that $U_{q,\theta+n-q}$ is Beta with parameters $(q, \theta + n - q)$, $\mu_\theta$ is a Gamma process with shape $\theta H$, $\{G_{j,n}\}$ are independent Gamma random variables with shape $e_{j,n}$ and scale 1. The distribution of $\mathbf{Y} = (\mathbf{Y}^*, \mathbf{p})$ is $\mathbb{P}(d\mathbf{Y}|\theta H)$. In particular, if $H$ is nonatomic, the $Y_j^*$ for $j = 1, \ldots, n(\mathbf{p})$ are i.i.d. $H$, and the distribution of $\mathbf{p}$ is $\pi(\mathbf{p}|\theta)$. Statement* (i) *implies the following results.*

(ii) *For all $\theta$ and $q = 1$,*

$$(19) \quad \mathcal{L}(\mu_{\theta,\theta-1}) = \mathcal{L}(U_{1,\theta}\mu_\theta + U_{1,\theta}T_1\delta_{Y_1}).$$

*If $\mu_\theta$ denotes a Gamma process with arbitrary shape parameter $\theta H$, then*

$$(20) \quad \mathcal{L}(\mu_\theta) = \mathcal{L}(U_{\theta,1}\mu_\theta + U_{\theta,1}T_1\delta_{Y_1}).$$

(iii) *For all positive $\theta$ and $q$,*

$$(21) \quad \mathcal{L}(\mu_{\theta,\theta-q}) = \mathcal{L}(T_q P_\theta),$$

*where $T_q$ is a Gamma random variable with shape $q$ and scale 1 independent of $P_\theta$, which is a Dirichlet process with shape $\theta H$. Hence for all positive $\theta$ and $q$, $\mathcal{L}(T_\theta \mu_{\theta,\theta-q}) = \mathcal{L}(T_q \mu_\theta)$, where $T_\theta$ is Gamma with shape $\theta$ and scale 1, independent of $\mu_{\theta,\theta-q}$. Similarly, $T_q$ and $\mu_\theta$ are independent.*

PROOF. The distributional identity in (i) is a direct consequence of the mixture representation of the law of $\mu_{\theta,\theta-q}$, in the form of the posterior distribution of $\mu_{\theta,\theta-q}|\mathbf{Y}$ and $\mathbb{P}(d\mathbf{Y}|\theta H)$, deduced from the expression for the Laplace functional in Theorem 2.1(ii). Note that all quantities on the right-hand side of (18), including $\mathbf{p}$, are random. We now show statement (iii) follows from statement (i). Notice $T_{\theta+n} := \mu_\theta(\mathcal{Y}) + \sum_{j=1}^{n(\mathbf{p})} G_{j,n}$ is a Gamma random variable with shape $\theta + n$ independent of $U_{q,\theta+n-q}$. Moreover, using the mixture representation of the Dirichlet process derived from its posterior distribution and $\mathbb{P}(d\mathbf{Y}|\theta H)$, it follows that $(\mu_\theta + \sum_{j=1}^{n(\mathbf{p})} G_{j,n}\delta_{Y_j^*})/T_{\theta+n}$ is a Dirichlet process with shape $\theta H$, independent of $T_{\theta+n}$ and $U_{q,\theta+n-q}$. Hence the right-hand side of (18) can be written as $U_{q,\theta+n-q}T_{\theta+n}P_\theta$. The result is completed by noting that $U_{q,\theta+n-q}T_{\theta+n}$ is equal in distribution to $T_q$. $\square$

REMARK 2.4. The distributional identities in (19) and (20), which are new, are analogous to similar identities for Dirichlet processes which have a variety of applications, as can be seen in Diaconis and Kemperman (1996), Sethuraman (1994) and Hjort (2003). In addition, distributional results between Beta and Gamma random variables have quite striking consequences, as can be seen from, for instance, Dufresne (1998). We view our results as functional extensions of some of those ideas, and it seems worthwhile to pursue more analogous results. Note importantly that our results do not require that $\mu_\theta(g)$ is Gamma distributed.

REMARK 2.5. The expression in Corollary 2.1 is obtained by evaluating the Laplace transform of the right-hand side of (19), in the order of integration of $\mu_\theta$, $T_1$, $Y_1$ and finally $U_{1,\theta}$. It is evident that other equivalent expressions can be formed by changing the order of integration. It is no coincidence that $U_{1,\theta}$ has the same distribution as $T_1/(T_1 + T)$ where $T = \int_{\mathcal{Y}} \mu(dy) = T_\theta$. Additionally, further representations can be obtained by using the distributional identity

$$U_{1,\theta} = \frac{(T_1)^p}{(T_1)^p + T_\theta \tau_p},$$

where $T_1, T_\theta$ and $\tau_p$ are all independent and $\tau_p$ is a stable random variable with index $0 < p < 1$.

2.1. *Distributional characterizations via the Beta-Gamma calculus.* The expression (21) tells us precisely that, for all $\theta$ and $q$, a Beta-Gamma process with parameters $\theta H$ and $\theta - q$ is equivalent in distribution to a Dirichlet process with shape $\theta H$, scaled by an independent Gamma random variable with shape $q$. Hence, using this interpretation the first result in Theorem 2.1 is an immediate consequence of

$$E[e^{-z\mu_{\theta,\theta-q}(g)}] = \frac{1}{\Gamma(q)} \int_0^\infty t^{q-1}\left[\int_{\mathcal{M}} e^{-ztP(g)}\mathcal{D}(dP|\theta H)\right]e^{-t}\,dt$$
$$= \int_{\mathcal{M}} (1 + zP(g))^{-q}\mathcal{D}(dP|\theta H).$$

Although this viewpoint at first may seem to have limited usage, it has interesting consequences when combined with our other results, within the context of the Beta-Gamma calculus. Note that since $\mu_{\theta,\theta-q}(g) = T_q P_\theta(g)$, one may apply the special features of the classical Beta-Gamma calculus combined with our results to deduce the following characterization of when $P_\theta(g)$ is Beta distributed, in the case where $g(Y)$ is not an indicator variable.

PROPOSITION 2.3. *Let $\theta$, $q$ and $\alpha$ denote positive real numbers. Suppose that for some $q > 0$, $\mu_{\theta,\theta-q}(g) = T_\alpha$, that is, it is Gamma distributed with shape parameter $\alpha$ and scale 1; then for the case $0 < \alpha < q$, $P_\theta(g) = U_{\alpha,q-\alpha}$ and $\mu_\theta(g) = T_\theta U_{\alpha,q-\alpha}$. From the distributional identity* (18), *this is true if one chooses the distribution of $g(Y)$ such that, for a fixed $n \geq 1$,*

$$(22) \qquad T_\theta U_{\alpha,q-\alpha} + \sum_{j=1}^{n(\mathbf{p})} G_j g(Y_j^*) = T_{\theta+n} U_{\alpha,q-\alpha},$$

*where $T_\theta U_{\alpha,q-\alpha}$ and $\sum_{j=1}^{n(\mathbf{p})} G_j g(Y_j^*)$ are independent. For clarity, when $n=1$,* (22) *specializes to $T_\theta U_{\alpha,q-\alpha} + T_1 g(Y_1) = T_{\theta+1} U_{\alpha,q-\alpha}$.*

Proposition 2.3 provides a characterization for the reverse question as to which choice of $H$ produces a Beta distribution for $P_\theta(g)$. This is of course seen to be equivalent to specifying $H$ to induce a particular distribution on the quantity $\sum_{j=1}^{n(\mathbf{p})} G_j g(Y_j^*)$, for some fixed value of $n$, which satisfies the constraint (22). It is clear by using (18), that this particular feature, of inducing a distribution on $\sum_{j=1}^{n(\mathbf{p})} G_j g(Y_j^*)$ satisfying appropriate constraints, can be applied to any choice of distribution for $P_\theta(g)$. However, because of the available independence properties between Beta and Gamma random variables, the occurrence of a Beta distribution for $P_\theta(g)$ can be checked several ways, not available to other distributions. In particular, note that $P_\theta(g) = U_{\alpha,q-\alpha}$ if and only if $\mu_{\theta,\theta-q} = T_\alpha$. Hence one can choose $g(Y)$ (or check this), such that the Laplace transform

$$E[e^{-z\mu_{\theta,\theta-q}(g)}] = (1+z)^{-\alpha},$$

or perhaps more easily using the Gamma process to check whether $\mu_\theta(g)$ satisfies

$$e^{-\int_{\mathcal{Y}} \log[1+zg(y)]\theta H(dy)} = E[e^{-zT_\theta U_{\alpha,q-\alpha}}].$$

We show in the next proposition the limitations, within the context of Proposition 2.3, of choosing $\sum_{j=1}^{n(\mathbf{p})} G_j g(Y_j^*)$ to be a Gamma random variable.

PROPOSITION 2.4. *Suppose that $H$ is chosen such that for some fixed $n \geq 1$, $\sum_{j=1}^{n(\mathbf{p})} G_j g(Y_j^*)$ has a Gamma distribution with shape parameter $c_n$, depending on $n$, and scale 1. Then the only value of $c_n$ such that $P_\theta(g)$ has a*

*Beta distribution is $c_n = n/2$. Moreover, the distribution of $P_\theta(g)$ is a symmetric Beta distribution with parameters $(\theta + n/2, \theta + n/2)$, for all $n \geq 1$ and all $\theta > 0$. Equivalently $\mu_{\theta,-(\theta+n)}(g) = T_{\theta+n/2}$ and $\mu_\theta(g) = T_\theta U_{\theta+n/2,\theta+n/2}$, and hence is never Gamma distributed. These specifications correspond to the choice of $\alpha = \theta + n/2$ and $q = 2\theta + n$ in Proposition* 2.3.

PROOF. The proof is obtained by applying Theorem 2 of Dufresne (1998), combined with the constraints deduced from (22) in Proposition 2.3. It follows that the only solution is given by

$$T_\theta U_{\theta+n/2,\theta+n/2} + T_{n/2} = T_{\theta+n} U_{\theta+n/2,\theta+n/2}. \qquad \square$$

It is already known [see Cifarelli and Melilli (2000)] that if $g(Y) = U_{1/2,1/2}$, the arcsine law, then $P_\theta(g) = U_{\theta+1/2,\theta+1/2}$. That is the case of $n=1$ in Proposition 2.4. The case for $n=2$ corresponds to

$$T_\theta U_{\theta+1,\theta+1} + [pT'_1 + T_1]g(Y_1) + (1-p)T'_1 g(Y_2) = T_{\theta+2} U_{\theta+1,\theta+1},$$

where $T'_1$, $T_1$ are independent exponential (1) random variables. $Y_1$ and $Y_2$ both have distribution $H$, but may be tied. $p$ is a Bernoulli random variable with success probability $1/(\theta+1)$, corresponding to the case where $Y_1 = Y_2$, from the Blackwell–MacQueen urn scheme. It is not immediately clear how to choose $H$ such that $[pT'_1 + T_1]g(Y_1) + (1-p)T'_1 g(Y_2) = T_1$.

**Acknowledgments.** I would like to thank Antonio Lijoi and Igor Prünster for re-stimulating my interest in this problem and elucidating some key points. Thanks to Jim Pitman for pointing out the Dufresne reference. Finally I would like to thank John W. Lau for his help in proofreading this manuscript.

Department of Information
and Systems Management
Hong Kong University of Science
and Technology
Clear Water Bay, Kowloon
Hong Kong
e-mail: lancelot@ust.hk